\documentclass[a4paper,12pt]{article}
\usepackage{geometry} 
\usepackage{amstext, graphicx, latexsym, amssymb, amscd, amsfonts, float,placeins,multirow,color,caption}
\usepackage{graphics}
\usepackage{amsthm}
\usepackage{amsfonts}
\usepackage{amsmath}\allowdisplaybreaks[3]
\usepackage[hidelinks]{hyperref}
\usepackage{longtable}
\usepackage{subcaption} 
\usepackage{mathtools}
\usepackage{multirow,lipsum}
\DeclareMathAlphabet{\mathpzc}{OT1}{pzc}{m}{it}
\usepackage[titletoc,toc,page]{appendix}
\theoremstyle{definition}
\renewcommand{\theequation}{\thesection.\arabic{equation}}
\newtheorem{theorem}{\bf Theorem}[section]

\newtheorem{example}{\bf Example}[section]
\numberwithin{equation}{section}
\numberwithin{table}{section}
\numberwithin{figure}{section}

\renewcommand\appendix{\par
\setcounter{section}{0}%
\setcounter{subsection}{0}%
\setcounter{table}{0}
\setcounter{figure}{0}
\setcounter{equation}{0}
\setcounter{theorem}{0}
\gdef\thetable{\Alph{table}}
\gdef\thefigure{\Alph{figure}}
\gdef\theequation{\Alph{section}-\arabic{equation}}
\section*{Appendix}
\gdef\thesection{\Alph{section}}
\setcounter{section}{1}}

\title{On risk models with dependence}
\author{Marjan Qazvini\footnote{School of Mathematical and Computer Sciences, Heriot-Watt University Malaysia, No 1 Jalan Venna P5/2, Precinct 5, 62200 Putrajaya, Malaysia, E-mail M.Qazvini@hw.ac.uk}\\\\
\small }

\date{} 

\begin{document}
\maketitle

\begin{abstract}

In this paper we consider the classical and Erlang(2) risk processes when the inter-claim times and claim amounts are dependent. We assume that the dependence structure is defined through a Farlie-Gumbel-Morgenstern (FGM) copula and show that the methods used to derive results in the classical risk model can be modified to derive results in a dependent risk process. We find expressions for the survival probability and the probability of maximum surplus before ruin.  

\end{abstract}

{\bf Keywords}: Classical risk model; Erlang(2) risk process; FGM copula; survival probability; maximum surplus before ruin

\section{Introduction}
The insurance operation is subject to the cyclical claim experience, meaning that the claim frequency is likely to rise in different seasons. For example, we expect that in winters more accidents to be reported to the insurance companies. Therefore, unlike the classical risk model which is based on the assumption that the inter-claim times and the amount of claims are independent, in real life we can find examples that these two random variables are dependent. Allowing for such dependence may lead to an increase in the probability of ruin.

\par There is much research that considers a structure that shows claim amounts and inter-claims are dependent. Beard et al. (1990) referred to this issue and provided different examples to justify such dependence. Albrecher and Boxma (2004) consider a situation that the distribution of the time between two claims depends on the previous claim size. Albrecher and Teugels (2006) study this setting by considering a copula for the dependence structure. They derive asymptotic results for the infinite and finite time ruin probabilities. Boudreault et al. (2006) consider a dependence structure such that the distribution of the next claim amount is defined in terms of the time elapsed since the last claim and find an explicit expression for the Laplace transform of the time of ruin. Cossette et al. (2010) extend their results to the situation that the claim amounts and the inter-claim times are dependent through a Farlie-Gumbel-Morgenstern (FGM) copula. Such dependence structure is also considered when the underlying risk process is Sparre Andersen. See, for example, Willmot and Woo (2012), Chadjiconstantinidis and Vrontos (2013), and Cossette et al. (2015). Jiang and Yang (2016) obtain an expression for the maximum surplus prior to ruin when the dependence structure is defined through a FGM copula.
\par In this paper we assume that the dependence structure between the inter-claim times and claim amounts is captured through a FGM copula. We remark that this copula does not show strong dependence between two random variables. However, It is widely used in the literature due to its simple functional form. Our purpose, in this study, is to show that the methods applied in the classical risk model can be adapted to derive results in models with dependence. In particular, we show that the methodology used by Dickson (1998) can be adapted to find the survival probability in the classical risk model considered by Cossette et al. (2010), the survival probability in an Erlang risk process considered by Chadjiconstantinidis and Vrontos (2013) and the probability of the maximum surplus before ruin studied by Jiang and Yang (2016). 
This paper is organised as follows. In the next section, we introduce models and notation. Section 3 considers the survival probability in the classical risk model and Section 4 considers the same problem in an Erlang(2) risk process. In Section 5 we find an expression for the probability of maximum surplus before ruin.

\section{Model and notation}\label{Sec:Model}
The risk surplus process in the classical risk model is defined by $U(t) = u+ct-S(t)$, where $u$ is the initial surplus, $S(t) = \sum_{i=1}^{N(t)} X_i$ is a compound Poisson process and $X_i$ is the amount of the $i$th claim with distribution function $F_X(x)$ and density $f_X(x)$ and $c$ is the premium rate. Also, $N(t)$ is a Poisson process. Let $\{W_i\}_{i=1}^{\infty}$ be a sequence of independent and identically distributed random variables and represent the inter-claim times. We assume $W$ is exponentially distributed with mean $1/\lambda$ and denote the distribution function of $W$ by $F_W(t)$ and its density by $f_W(t)$. To ensure that ruin will not occur almost surely, throughout we assume that the positive loading condition i.e. $E[c W_i - X_i] > 0$ for $i=1,2,\dots$ satisfies. Let $T$ denote the time of ruin and is defined as $T = \inf\{t: U(t) < 0 | U(0) = u\}$ with $T=\infty$ if $U(t)\geq 0$ for all $t>0$. Then the ultimate ruin probability is defined as $\psi(u) = \textnormal{Pr}(T<\infty | U(0) =u) = 1-\phi(u)$, where $\phi(u)$ is the survival probability in the classical risk model.

\par In the classical risk model we assume that $X_i$ and $W_i$ are independent. Here, we relax this assumption and consider a sequence of independent and identically distributed random vectors $\{(X_i, W_i)\}_{i=1}^{\infty}$ with joint density $f_{X,W}(x,t)$ and joint distribution function $F_{X,W}(x,t)$. Throughout we denote the Laplace transform of a function  $\alpha(u)$ by $\tilde{\alpha}(s)$.
\par To build the joint distribution of $(X,W)$ in terms of a copula, we need to apply the Sklar's theorem. (See, for example, Nelsen, 2006.)
\begin{theorem}[\bf Sklar's theorem]
If $(X,W)$ form a continuous random vector with continuous distribution $F_X$ and $F_W$ and joint distribution function $F_{X,W}$, then there exists a unique copula $C_{\theta}$ which is defined on $[0,1]^2$ with uniform margins such that
\begin{eqnarray}
F_{X,W}(x,t) = C_{\theta}(F_X(x), F_W(t))\nonumber
\end{eqnarray}
with density 
\begin{eqnarray}
f_{X,W}(x,t) = c_{\theta}(F_X(x),F_W(t)) f_X(x) f_W(t)\nonumber 
\end{eqnarray}
where $\theta$ is the copula's parameter and $c_{\theta}(F_X(x), F_W(t)) = \frac{\partial^2}{\partial F_X \partial F_W}C(F_X(x),F_W(t))$.
\end{theorem}
Similar to Cossette et al. (2010) we define the joint distribution of $(X,W)$ in terms of a FGM copula, which is a bivariate, one-parameter copula and is given by

\begin{eqnarray}
C_{\theta}(u,v) = uv +  \theta u v (1-u) (1-v)\nonumber
\end{eqnarray}
with the density being
\begin{eqnarray}
c_{\theta}(u,v) = 1 + \theta (1-2 u) (1-2v)\nonumber.
\end{eqnarray}
For the properties of this copula, readers are referred to Nelsen (2006).
\par Noting that $u = F_X$ and $v = F_W$, we can write the joint density of $X$ and $W$ as
\begin{eqnarray}\label{eq:jointWX}
f_{X,W}(x,t) = \left[1 + \theta (1-2 F_X(x))(1-2F_W(t))\right] f_X(x) f_W(t).
\end{eqnarray}
Substituting for $F_W$ and $f_W$ we obtain
\begin{eqnarray}\label{joint-X-W}
f_{X,W}(x,t) &=& \left[1 + \theta (1-2F_X(x)) (2 e^{-\lambda}-1)\right] f_X(x) \lambda e^{-\lambda t}\nonumber\\
&=& f_X(x) \lambda e^{-\lambda t} + \theta f_X(x) \lambda e^{-\lambda t} (1-2 F_X(x)) (2 e^{-\lambda t}-1)\nonumber\\
&=& f_X(x) \lambda e^{-\lambda t} + \theta \lambda e^{-\lambda t} h_X(x) (2 e^{-\lambda t}-1),
\end{eqnarray}
where $h_X(x) = f_X(x) (1-2F_X(x))$. 

\section{Survival probability}
In this section, we apply the method in Dickson (1998) to derive the Laplace transform of the survival probability. Then we find numerical expressions for the survival probability by inverting its Laplace transform. 
\par Conditioning on the time and the amount of the first claim, we have
\begin{eqnarray}
\phi(u)&=& \int_0^{\infty} \int_0^{u+ct} f_{X,W}(x,t) \phi(u+ct-x) dx dt\nonumber\\
&=&\int_0^{\infty} \lambda e^{-\lambda t} \int_0^{u+ct}f_X(x)\phi(u+ct-x)dx dt\nonumber\\
&&+ \theta \int_0^{\infty} \lambda e^{-\lambda t} (2e^{-\lambda t}-1) \int_0^{u+ct}h_X(x)\phi(u+ct-x)dx dt\nonumber\\
&=&\int_0^{\infty}\lambda e^{-\lambda t}\int_0^{u+ct}f_X(x)\phi(u+ct-x) dx dt+ 2 \theta \int_0^{\infty} \lambda e^{-2\lambda t} \int_0^{u+ct} h_X(x)\phi(u+ct-x)dx dt\nonumber\\
&&-\theta \int_0^{\infty} \lambda e^{-\lambda t}\int_0^{u+ct}h_X(x)\phi(u+ct-x)dx dt\nonumber.
\end{eqnarray}
Substituting $s=u+ct$, gives
\begin{eqnarray}
\phi(u)&=&\frac{\lambda}{c}\int_u^{\infty} e^{-\lambda(s-u)/c} \int_0^s f_X(x)\phi(s-x)dx ds+\frac{2\theta \lambda}{c}\int_u^{\infty} e^{-2\lambda(s-u)/c}\int_0^s h_X(x)\phi(s-x) dx ds\nonumber\\
&&-\frac{\theta \lambda}{c}\int_u^{\infty} e^{-\lambda(s-u)/c} \int_0^s h_X(x)\phi(s-x)dx ds\nonumber.
\end{eqnarray}
Differentiating with respect to $u$, we have
\begin{eqnarray}
\frac{d}{du}\phi(u)&=& \frac{\lambda^2}{c^2} \int_u^{\infty} e^{-\lambda(s-u)/c}\int_0^s f_X(x)\phi(s-x)dx ds -\frac{\lambda}{c}\int_0^u f_X(x)\phi(u-x)dx\nonumber\\
&&+\frac{4\theta \lambda^2}{c^2}\int_u^{\infty}e^{-2\lambda(s-u)/c}\int_0^s h_X(x)\phi(s-x)dx ds-\frac{2\theta \lambda}{c}\int_0^u h_X(x)\phi(u-x)dx\nonumber\\
&&-\frac{\theta \lambda^2}{c^2}\int_u^{\infty}e^{-\lambda(s-u)/c}\int_0^s h_X(x)\phi(s-x)dx ds+\frac{\theta \lambda}{c}\int_0^u h_X(x)\phi(u-x)dx\nonumber\\
&=&\frac{\lambda}{c} \phi(u)+\frac{2\theta \lambda^2}{c^2}\int_u^{\infty} e^{-2\lambda(s-u)/c}\int_0^s h_X(x)\phi(s-x)dx ds-\frac{\lambda}{c}\int_0^u f_X(x)\phi(u-x)dx\nonumber\\
&&-\frac{\theta \lambda}{c}\int_0^u h_X(x)\phi(u-x)dx\nonumber\\
&=&\frac{\lambda}{c} \phi(u)+\frac{2\theta \lambda^2}{c^2}\int_u^{\infty} e^{-2\lambda(s-u)/c}\int_0^s h_X(x)\phi(s-x)dx ds-\frac{\lambda}{c}A(u)-\frac{\theta \lambda}{c}B(u)\nonumber,
\end{eqnarray}
where $A(u) = \int_0^u f_X(x)\phi(u-x)dx$ with the Laplace transform $\tilde{A}(s)=\tilde{f}(s)\tilde{\phi}(s)$ and $B(u)= \int_0^u h_X(x)\phi(u-x)dx$ with the Laplace transform $\tilde{B}(s)=\tilde{h}(s)\tilde{\phi}(s)$. Differentiating again, we get
\begin{eqnarray}\label{eq:DDS}
\frac{d^2}{du^2}\phi(u)&=& \frac{\lambda}{c}\frac{d}{du}\phi(u)+\frac{4\theta\lambda^3}{c^3}\int_u^{\infty}e^{-2\lambda(s-u)/c}\int_0^s h_X(x)\phi(s-x)dx ds-\frac{2\theta \lambda^2}{c^2}B(u)\nonumber\\
&&-\frac{\lambda}{c}\frac{d}{du}A(u)-\frac{\theta \lambda}{c}\frac{d}{du}B(u)\nonumber\\
&=&\frac{3\lambda}{c}\frac{d}{du}\phi(u)-\frac{2\lambda^2}{c^2}\phi(u)+\frac{2\lambda^2}{c^2}A(u)-\frac{\lambda}{c}\frac{d}{du}A(u)-\frac{\theta\lambda}{c}\frac{d}{du}B(u).
\end{eqnarray}
Taking the Laplace transform from (\ref{eq:DDS}) we obtain
\begin{eqnarray}
&&c^2(s^2 \tilde{\phi}(s)-s\phi(0)-\phi^{'}(0))-3\lambda c (s\tilde{\phi}(s)-\phi(0))+2\lambda^2\tilde{\phi}(s)\nonumber\\
&=&2\lambda^2\tilde{f}(s)\tilde{\phi}(s)-\lambda c (s\tilde{f}(s)\tilde{\phi}(s)-0)-\theta\lambda c(s\tilde{h}(s)\tilde{\phi}(s)-0)\nonumber.
\end{eqnarray}
So,
\begin{eqnarray}
\tilde{\phi}(s)=\frac{c^2s\phi(0)+c^2\phi^{'}(0)-3\lambda c\phi(0)}{c^2 s^2-3\lambda c s+2\lambda^2(1-\tilde{f}(s))+\lambda cs(\tilde{f}(s)+\theta \tilde{h}(s))}\nonumber.
\end{eqnarray}
Define $L$ to be the maximum aggregate loss, then we have $\textnormal{Pr}(L \leq u) = \phi(u)$ and $E[e^{-sL}]=s \tilde{\phi}(s)$ with $\lim_{s \rightarrow 0} E[e^{-sL}] = 1$. Using L'Hospital's rule we have $c^2 \phi^{'}(0)-3\lambda \phi(0) = -2\lambda c+2\lambda^2 m_1$, where $m_1$ is the mean of the claim amount distribution.  Therefore the Laplace transform of the survival probability is given by
\begin{eqnarray}\label{eq:LTS}
\tilde{\phi}(s) = \frac{c^2 s\phi(0) -2 \lambda c+2\lambda^2 m_1}{c^2 s^2-3\lambda c s+2\lambda^2(1-\tilde{f}(s))+\lambda cs(\tilde{f}(s)+\theta \tilde{h}(s))},
\end{eqnarray}
which can be easily inverted using software like Mathematica. 

\begin{example}
We consider the situation when both claim amounts and inter-claim times follow an exponential distribution with mean $1$ and the premium rate $1.5$. We find the survival probability for different values of $\theta$ by inverting its Laplace transform.
\begin{itemize}
\item $\theta = -0.5$:
\begin{eqnarray}
\phi(u) &=& 1- (0.5747 + 0.3848 \phi(0)) e^{-R_1u} + (0.0042+0.0200 \phi(0)) e^{-R_2 u} \nonumber\\
&&- (0.4295 - 1.3649 \phi(0)) e^{-R_3 u}\nonumber,
\end{eqnarray}
where $R_1 = 0.2976$, $R_2 = 2.1148$ and $R_3 = -1.4123$. Since Lundberg's inequality applies, the coefficient of $e^{-R_3 u}$ must be $0$. Hence $\phi(0) = 0.3147$. Similarly, we can find the survival probability for other values of $\theta$. 
\item $\theta=0$ and $\phi(0) = 0.3333$:
\begin{eqnarray}
\phi(u) = 1- 0.6667 e^{-0.3333 u}\nonumber,
\end{eqnarray}
\item $\theta =0.5$ and $\phi(0) = 0.3548$:
\begin{eqnarray}
\phi(u) = 1 - 0.6311 e^{-0.3788 u} - 0.0140 e^{-1.8736 u}\nonumber.
\end{eqnarray}
\end{itemize}
\end{example}
This example is based on the parameters used in Cossette et al. (2010). As we can see our method can be simply applied to find the survival probability.

\section{Survival probability with Erlang$(2,2)$ inter-arrival times}
In this section, we consider a risk process where $W_i$ is distributed according to an Erlang$(2,2)$ distribution with density function $f_W(t) = \beta^2 t e^{-\beta t}$. We denote the survival probability in this model by $\delta(u)$. First we find an expression for $f_{X,W}$ by applying (\ref{eq:jointWX}). Hence
\begin{eqnarray}
f_{X,W}(x,t) = f_X(x) f_W(t) + \theta h_X(x) k(t)\nonumber,
\end{eqnarray}
where $h_X$ is as defined in Section \ref{Sec:Model} and $k_W(t) = f_W(t) (1-2F_W(t))$. We also use the fact that 
$f^{'}_W(t) = \beta^2 e^{-\beta t} - \beta f_W(t)$ and 
\begin{eqnarray}
k_W^{'}(t) &=& f^{'}_W(t) - 2 f^{'}_W(t) F_W(t) - 2 f^2_W(t)\nonumber\\
&=& f^{'}_W(t) (1-2F_W(t)) - 2 f^2_W(t)\nonumber\\
&=& (\beta^2 e^{-\beta t} - \beta f_W(t)) (1-2F_W(t)) - 2 f^{2}_W(t)\nonumber\\
&=& \beta^2 e^{-\beta t} (1-2F_W(t)) - \beta k_W(t) - 2 f^2_W(t)\nonumber.
\end{eqnarray}

\par We now condition on the time and the amount of the first claim and write
\begin{eqnarray}
\delta(u) &=& \int_0^{\infty} \int_0^{u+ct} f_{X,W}(x,t) \delta(u+ct-x)dx dt\nonumber\\
&=& \int_0^{\infty} f_W(t) \int_0^{u+ct} f_X(x) \delta(u+ct-x)  dx dt\nonumber\\
&&+\theta \int_0^{\infty} k_W(t) \int_0^{u+ct} h_X(x) \delta(u+ct-x) dx dt\nonumber.
\end{eqnarray}
Substituting $s=u+ct$ gives
\begin{eqnarray}
c \delta(u) &=& \int_u^{\infty} f_W\left(\frac{s-u}{c}\right) \int_0^s f_X(x) \delta(s-x) dx ds\nonumber\\
&& + \theta \int_u^{\infty} k_W\left(\frac{s-u}{c}\right) \int_0^s h_X(x) \delta(s-x) dxds.
\end{eqnarray}
Differentiating with respect to $u$, we have
\begin{eqnarray}
c \frac{d}{du} \delta(u) &=& \beta \delta(u) - \frac{1}{c} \int_u^{\infty} \beta^2 e^{-\beta(s-u)/c} \int_0^s f_X(x) \delta(s-x) dx ds\nonumber\\
&&- \frac{\theta}{c} \int_u^{\infty} \beta^2 e^{-\beta(s-u)/c} \int_0^s h_X(x) \delta(s-x) dx ds\nonumber\\
&&+\frac{2\theta}{c} \int_u^{\infty} \beta^2 e^{-\beta(s-u)/c} F_W\left(\frac{s-u}{c}\right) \int_0^s h_X(x) \delta(s-x) dx ds\nonumber\\
&&+\frac{2\theta}{c} \int_u^{\infty} f^2_W\left(\frac{s-u}{c}\right) \int_0^s h_X(x) \delta(s-x) dx ds\nonumber
\end{eqnarray}
and differentiating again, yields
\begin{eqnarray}
c^2 \frac{d^2}{du^2}\delta(u)&=& 2 \beta c \frac{d}{du} \delta(u) - \beta^2 \delta(u) + \beta^2 \int_0^u f_X(x) \delta(u-x) dx+ \theta \beta^2 \int_0^u h_X(x) \delta(u-x) dx\nonumber\\
&&- \frac{2\theta \beta}{c} \int_u^{\infty} f^2_W\left(\frac{s-u}{c}\right) \int_0^s h_X(x) \delta(s-x) dx ds\nonumber\\
&&-\frac{2\theta}{c} \int_u^{\infty} \beta^2 e^{-\beta(s-u)/c} f_W\left(\frac{s-u}{c}\right) \int_0^s h_X(x) \delta(s-x) dx ds\nonumber\\
&&-\frac{4\theta}{c} \int_u^{\infty} f_W\left(\frac{s-u}{c}\right) f^{'}_W\left(\frac{s-u}{c}\right) \int_0^s h_X(x) \delta(s-x) dx ds,\nonumber
\end{eqnarray}
which after inserting for $f_W$ can also be written as
\begin{eqnarray}\label{eq:secondd}
c^2 \frac{d^2}{du^2} \delta(u) - 2 \beta c \frac{d}{du}\delta(u) + \beta^2 \delta(u) &=& \beta^2 \int_0^u f_X(x) \delta(u-x) dx+ \theta \beta^2 \int_0^u h_X(x) \delta(u-x) dx\nonumber\\
&&+ \frac{2 \theta \beta^5}{c^3} \int_u^{\infty} e^{-2\beta(s-u)/c} (s-u)^2 \eta(s) ds\nonumber\\
&&-\frac{6 \theta \beta^4}{c^2} \int_u^{\infty} e^{-2\beta(s-u)/c} (s-u) \eta(s) ds,
\end{eqnarray}
where $\eta(s) = \int_0^s h_X(x) \delta(s-x) dx$. We now recall the Dickson-Hipp operator and some of its properties. (See Dickson and Hipp, 2001, and Li and Garrido, 2004).

\par The operator $T_r $ for an integrable function $f$ and real $r$ is defined by
\begin{eqnarray}
T_r f(x) = \int_x^{\infty} e^{-r(u-x)} f(u) du.\nonumber
\end{eqnarray}
The operator $T_r$ satisfies the following properties:
\begin{itemize}
\item[(i)] $ T_r f(0) = \int_0^{\infty} e^{-r u} f(u) du = \tilde{f}(r)$. 
\item[(ii)] $T_r^n f(x) = \int_x^{\infty} e^{-r(u-x)} \frac{(u-x)^{n-1}}{(n-1)!} f(u) du$.
\item[(iii)] $T_s T^n_r f(0) = \frac{\tilde{f}(s)}{(r-s)^n} - \sum_{j=1}^n \frac{T^j_r f(0)}{(r-s)^{n+1-j}}$.
\end{itemize}
 
We can now rewrite (\ref{eq:secondd}) in terms of the Dickson-Hipp operator as

\begin{eqnarray}\label{eq:second1}
c^2 \frac{d^2}{d u^2} \delta(u) - 2 \beta c \frac{d}{du} \delta(u) + \beta^2 \delta(u) &=& \beta^2 \int_0^u f_X(x) \delta(u-x) dx + \theta \beta^2 \int_0^u h_X(x) \delta(u-x) dx\nonumber\\
&&+ \frac{2 \theta \beta^5}{c^3} 2 T^3_{2\beta/c} \eta(u) - \frac{6 \theta \beta^4}{c^2} T^2_{2\beta/c} \eta(u).
\end{eqnarray} 
Applying property (ii) we can take the Laplace transform from both sides of (\ref{eq:second1}), giving
\begin{eqnarray}\label{eq:second2}
&&c^2 (s^2 \tilde{\delta}(s) -s \delta(0) - \delta^{'}(0)) - 2 \beta c (s \tilde{\delta}(s) - \delta(0)) + \beta^2 \tilde{\delta}(s)\nonumber\\
&=& \beta^2 \tilde{f}(s) \tilde{\delta}(s) + \theta \beta^2 \tilde{h}(s) \tilde{\delta}(s) + \frac{2\theta \beta^5}{c^3} 2 T_s T^3_{2\beta/c} \eta(0) - \frac{6 \theta \beta^4}{c^2} T_s T^2_{2\beta/c} \eta(0).
\end{eqnarray}
Applying property (iii) and noting that $\eta(0)=0$ and that $\tilde{\eta}(s) = \tilde{h}(s) \tilde{\delta}(s)$, the Laplace transform of the survival property is given by
\begin{eqnarray}\label{eq:SPaL}
\tilde{\delta}(s) = \frac{c^2 s \delta(0) + c^2 \delta^{'}(0) -2 \beta c \delta(0)}{c^2 s^2 - 2 \beta c s + \beta^2 (1-\tilde{f}(s))- \theta \beta^2 \tilde{h}(s) -4 \theta \beta^5 \tilde{h}(s)  /(2\beta-cs)^3 + 6 \theta \beta^4 \tilde{h}(s) /(2\beta-cs)^2 }\nonumber\\
\end{eqnarray}
where the denominator is the Lundberg's equation given by formula (11) in Chadjiconstantinidis and Vrontos (2013) for $\delta=0$ and $n=2$. They show that for $\delta=0$ and $\theta \neq 0$, the Lundberg's equation has 
$3n-1$ roots $\rho_1, \dots, \rho_{3n-1}$ with positive real part including one equal to $0$ and two roots $-R_i$ with $Re(R_i) >0$, for $i=1,2$.

To eliminate $\delta^{'}(0)$ in numerator, we can apply the final value theorem; see Dickson (1998). According to this theorem, we have 
$\lim_{u \rightarrow \infty} \delta(u) = \lim_{s \rightarrow 0} s \tilde{\delta}(s)$, which in our case is interpreted as $\lim_{s \rightarrow 0} s \tilde{\delta}(s) = 1$. Hence

\begin{eqnarray}
s \tilde{\delta}(s) = \frac{c^2 s^2 \delta(0) + c^2 s \delta^{'}(0) -2 \beta c s\delta(0)}{c^2 s^2 -2 \beta c s + \beta^2 (1-\tilde{f}(s)) - \theta \tilde{h}(s) (\beta^2 + 4 \beta^5 /(2\beta -c s)^3 - 6 \beta^4/(2\beta -c s)^2)}\nonumber.
\end{eqnarray}

Applying L'Hospital's rule, we can find the limit as $s \rightarrow 0$. Therefore
\begin{eqnarray}
c^2 \delta^{'}(0) - 2 \beta c \delta(0) = -2 \beta c + \beta^2 m_1\nonumber
\end{eqnarray} 
and so
\begin{eqnarray}\label{eq:LTSpar}
\tilde{\delta}(s) = \frac{c^2 s \delta(0) -2 \beta c + \beta^2 m_1}{c^2 s^2 - 2\beta c s + \beta^2 (1-\tilde{f}(s)) - \theta \tilde{h}(s) (\beta^2 + 4\beta^5 /(2\beta-c s)^3 + 6 \beta^4/(2\beta -cs)^2)}.
\end{eqnarray}

In the next section we provide numerical examples and find expressions for $\delta(u)$.  

\begin{example}
We consider the situation when claim amounts follow an exponential distribution with mean $1$, $c=1.5$ and $\beta=2$. Inserting our parameter values in (\ref{eq:LTSpar}) and inverting, for different values of $\theta$ we find:

\begin{itemize}
\item $\theta = -1$:
\begin{eqnarray}
&&\delta(u)\nonumber\\
&=& 1 - (0.5637+0.2211 \delta(0)) e^{-0.3488 u} + (0.0090+0.0219 \delta(0)) e^{-2.1517 u}\nonumber\\
&&+(-0.0481+0.1975 \delta(0)) e^{3.6476 u} + (-0.1372+0.2779 \delta(0)) e^{1.8011 u}\nonumber\\
&&+((-0.1300-0.0133i)+(0.3620-0.1295i) \delta(0)) e^{(2.3592-1.1277i)u}\nonumber\\
&&+((-0.1300+0.0133i)+(0.3620+0.1295i)\delta(0)) e^{(2.3593+1.1277i)u}\nonumber.
\end{eqnarray}
To find $\delta(0)$ we note that the roots of the Lundberg's equation are $\rho_{1} = 3.6476, \rho_2 = 2.3592 + 1.1277i, \rho_3 = 2.3593-1.1277i, \rho_4 = 0$ and $-R_1 =  2.1517, -R_2 = 1.8011$. Therefore, we need to set the coefficients of $e^{\rho_i u}$, $i=1,\dots 3$ equal to zero. If we let $u \rightarrow \infty$, we have $\delta(u) =1$ and $e^{\rho_i u} \rightarrow \infty$. Therefore all $e^{\rho_i u}$ are very large number. So, equating all coefficients of $e^{\rho_i u}$ collectively to zero we find $\delta(0) = 0.3713$. So
\begin{eqnarray}
\delta(u) = 1 - 0.6458 e^{-0.3488 u} + 0.0171 e^{-2.1517 u}\nonumber.
\end{eqnarray}
Similarly, we can find the survival probability for other values of $\theta$ as
\item $\theta =-0.5$ and $\delta(0) = 0.3963$:
\begin{eqnarray} 
\delta(u) = 1 - 0.6134 e^{-0.3833 u} + 0.0098 e^{-2.0792 u}\nonumber,
\end{eqnarray}
\item $\theta = 0.5$ and $\delta(0) = 0.4579$:
\begin{eqnarray}
\delta(u) = 1 - 0.5289 e^{-0.4762 u} - 0.0132 e^{-1.9119 u}\nonumber,
\end{eqnarray}
\item $\theta =1$ and $\delta(0) = 0.4957$:
\begin{eqnarray}
\delta(u) = 1 - 0.4723 e^{-0.5410 u} - 0.0320 e^{-1.8116 u}\nonumber.
\end{eqnarray}
\end{itemize}
\end{example}

\section{The maximum surplus before ruin}
For $0 \leq u < b$, define
\begin{eqnarray}
\xi(u,b) = \textnormal{P}\left(\sup_{0 \leq t \leq T} U(t) < b, T< \infty | U(0) = u\right)\nonumber
\end{eqnarray}
to be the probability that ruin occurs from initial surplus $u$ without the surplus process reaching level $b$ prior to ruin. Also, let $\chi(u,b)$ denote the probability that the surplus process reaches the level $b>u$ from initial surplus $u$ without first falling below zero. Since eventually either ruin occurs without the surplus process attaining level $b$ or the surplus attains level $b$, we have $\chi(u,b) = 1-\xi(u,b)$. In this section we find a numerical expression for $\chi(u,b)$ by solving a differential equation.  
\par Conditioning on the time and the amount of the first claim, we have
\begin{eqnarray}
\chi(u,b)= \int_0^{\tau} \int_0^{u+ct} f_{X,W}(x,t) \chi(u+ct-x,b) dx dt+\int_{\tau}^{\infty} \int_0^{\infty} f_{X,W}(x,t) dx dt\nonumber,
\end{eqnarray}
where $u+c\tau = b$, i.e. $\tau$ is the time that the surplus process reaches $b$ if no claim occurs by time $\tau$. Inserting for $f_{X,W}$ and setting $s=u+ct$ gives
\begin{eqnarray}\label{eq:chid}
\chi(u,b)&=& \frac{\lambda}{c}\int_u^b e^{-\lambda(s-u)/c} \int_0^s f_X(x) \chi(s-x,b) dx ds+\frac{2\theta \lambda}{c}\int_u^b e^{-2\lambda(s-u)/c}\int_0^s h_X(x) \chi(s-x,b) dx ds\nonumber\\
&&-\frac{\theta \lambda}{c}\int_u^b e^{-\lambda(s-u)/c} \int_0^s h_X(x) \chi(s-x,b)dx ds+\frac{\lambda}{c}\int_b^{\infty} e^{-\lambda(s-u)/c}ds.
\end{eqnarray}
Differentiating twice with respect to $u$ yields
\begin{eqnarray}\label{eq:chi1}
\frac{d^2}{d u^2}\chi(u,b)= \frac{3\lambda}{c}\frac{d}{du}\chi(u,b)-\frac{2\lambda^2}{c^2}\chi(u,b)+\frac{2\lambda^2}{c^2}\mathpzc{A}(u)-\frac{\lambda}{c}\frac{d}{du}\mathpzc{A}(u)-\frac{\theta \lambda}{c}\frac{d}{du}\mathpzc{B}(u),
\end{eqnarray}
where $\mathpzc{A}(u) = \int_0^u f_X(x) \chi(u-x,b)dx$ and $\mathpzc{B}(u) = \int_0^u h_X(x) \chi(u-x,b)dx$. Assuming $f_X(x) = \alpha e^{-\alpha x}$ and hence $h_X(x) = \alpha e^{-\alpha x}(2 e^{-\alpha x}-1)$, we now substitute for $f_X$ and $h_X$ in $\mathpzc{A}$ and $\mathpzc{B}$ and differentiate (\ref{eq:chi1}) another two times to eliminate the integral terms, which gives
\begin{eqnarray}
\frac{d^4}{d u^4}\chi(u,b)&=& \left(\frac{3\lambda}{c}-3\alpha\right)\frac{d^3}{d u^3}\chi(u,b)+\left(\frac{8\alpha \lambda}{c}-2 \alpha^2 - \frac{2 \lambda^2}{c^2}-\frac{\alpha \lambda \theta}{c}\right)\frac{d^2}{d u^2}\chi(u,b)\nonumber\\
&&+\left(\frac{4 \alpha^2 \lambda}{c}-\frac{4\alpha \lambda^2}{c^2}\right)\frac{d}{du}\chi(u,b)
\end{eqnarray}
with the characteristic equation
\begin{eqnarray}
s^4 - \left(\frac{3\lambda}{c}-3\alpha\right) s^3 - \left(\frac{8\alpha \lambda}{c}-2\alpha^2-\frac{2\lambda^2}{c^2}-\frac{\alpha \lambda \theta}{c}\right) s^2 - \left(\frac{4\alpha^2 \lambda}{c}-\frac{4 \alpha \lambda^2}{c^2}\right) s =0\nonumber
\end{eqnarray}
and general solution being given by
\begin{eqnarray}
\chi(u,b) = a_0 + \sum_{i=1}^3 a_i e^{-R_i u},
\end{eqnarray}
where $\{R_i\}_{i=1}^{3}$ are the solutions to the characteristic equation. To find $\{a_i\}_{i=0}^3$ we need four equations. The first one is our boundary condition $\chi(b,b) =1$, giving
\begin{eqnarray}
1 = a_0 + \sum_{i=1}^3 a_i e^{-R_i b}.\nonumber
\end{eqnarray}
To obtain another equation, we substitute in the derivative of (\ref{eq:chid}) which is
\begin{eqnarray}\label{eq:chif}
\frac{d}{du}\chi(u,b)-\frac{\lambda}{c}\chi(u,b)&=& \frac{2\theta \lambda^2}{c^2}\int_u^b e^{-2\lambda(s-u)/c} \int_0^s h_X(x)\chi(s-x,b) dx ds\nonumber\\
&& - \frac{\lambda}{c}\int_0^u f_X(x) \chi(u-x,b) dx -\frac{\theta \lambda}{c}\int_0^u h_X(x) \chi(u-x,b) dx.\nonumber\\
\end{eqnarray}
Inserting for $\chi(u,b)$ in the left-hand side of (\ref{eq:chif}), we get
\begin{eqnarray}\label{eq:Lefth}
-\frac{\lambda}{c}a_0 - \sum_{i=1}^3 \left(\frac{\lambda}{c}+R_i\right) a_i e^{-R_i u}.
\end{eqnarray}
Inserting for $\chi(u,b)$, $f_X$ and $h_X$ in the right-hand side of (\ref{eq:chif}), the first term becomes
\begin{eqnarray}
&&\frac{a_0 c}{\alpha c+2 \lambda}\left(e^{-\alpha u} - e^{-\alpha b} e^{-2\lambda \tau}\right)-\frac{a_0 c}{2\alpha c+2\lambda}\left(e^{-2\alpha u}-e^{-2\alpha b}e^{-2\lambda \tau}\right)\nonumber\\
&&+2\alpha c\sum_{i=1}^3 \frac{ a_i \left(e^{-R_i u}- e^{-R_i b}e^{-2\lambda \tau}\right)}{(2\alpha-R_i)(R_i c+2\lambda)}  - 2\alpha c \sum_{i=1}^3 \frac{a_i \left(e^{-2\alpha u}-e^{-2\alpha b} e^{-2\lambda \tau}\right)}{(2\alpha-R_i)(2\alpha c+2\lambda)}\nonumber\\
&&-\alpha c\sum_{i=1}^3 \frac{a_i \left(e^{-R_i u}-e^{-R_i b}e^{-2\lambda \tau}\right)}{(\alpha-R_i)(R_i c+2\lambda)}+ \alpha c \sum_{i=1}^3 \frac{a_i \left(e^{-\alpha u}-e^{-\alpha b}e^{-2\lambda \tau}\right)}{(\alpha-R_i)(\alpha c+2\lambda)}.
\end{eqnarray}
We do the same thing for the other two integrals and respectively obtain
\begin{eqnarray}
a_0 - a_0 e^{-\alpha u} + \alpha \sum_{i=1}^3 \frac{a_i}{\alpha-R_i} \left(e^{-R_i u}-e^{-\alpha u}\right)
\end{eqnarray}
and
\begin{eqnarray}
a_0 \left(e^{-\alpha u} - e^{-2\alpha u}\right) + \sum_{i=1}^3 \frac{2 \alpha a_i}{2\alpha-R_i}\left(e^{-R_i u}-e^{-2\alpha u}\right)- \sum_{i=1}^3 \frac{\alpha a_i}{\alpha-R_i}\left(e^{-R_i u}-e^{-\alpha u}\right).
\end{eqnarray}
We note that (\ref{eq:Lefth}) only contains the terms $e^{-R_i}$. Therefore we need to eliminate the terms involving $e^{-\alpha u}$, $e^{-2 \alpha u}$ and $e^{-2 \lambda \tau}$ on the right-hand side of (\ref{eq:chif}). Hence we have another three equations by equating coefficients of $e^{-\alpha u}$ to zero, i.e.
\begin{eqnarray}
2 \theta \lambda \left(\frac{\alpha}{\alpha c+2\lambda} + \sum_{i=1}^3 \frac{\alpha c a_i}{(\alpha-R_i)(\alpha c+2 \lambda)}\right) + \left(a_0 + \sum_{i=1}^3 \frac{\alpha a_i}{\alpha-R_i}\right)-\theta \left(a_0+\sum_{i=1}^3 \frac{\alpha a_i}{\alpha-R_i}\right)=0\nonumber,
\end{eqnarray}
by equating coefficients of $e^{-2\alpha u}$ to zero, i.e.
\begin{eqnarray}
-2 \lambda \left(\frac{\alpha}{2\alpha c+2\lambda} + \sum_{i=1}^3\frac{\alpha c a_i}{(2\alpha-R_i)(\alpha c+\lambda)}\right)+ \left(a_0 + \sum_{i=1}^3 \frac{2\alpha a_i}{2\alpha-R_i}\right) = 0,\nonumber
\end{eqnarray}
and by equating coefficients of $e^{-2\lambda \tau}$ to zero, i.e
\begin{eqnarray}
&&\frac{-a_0 e^{-\alpha b}}{\alpha c+2 \lambda} + \frac{a_0 e^{-2\alpha b}}{2\alpha +2 \lambda}- \sum_{i=1}^3 \frac{2\alpha a_i e^{-R_i b}}{(2\alpha-R_i)(R_i c+ 2\lambda)}+ \sum_{i=1}^3 \frac{\alpha a_i e^{-2 \alpha b}}{(2\alpha-R_i)(\alpha c+\lambda)}\nonumber\\
&&+\sum_{i=1}^3\frac{\alpha a_i e^{-R_i b}}{(\alpha-R_i)(R_ic+2\lambda)}-\frac{\alpha a_i e^{-\alpha b}}{(\alpha-R_i)(\alpha c+2\lambda)}= 0\nonumber.
\end{eqnarray}

\begin{example}
We consider the situation when both inter-claim times and claim amounts follow an exponential distribution with mean $1$ and the premium rate $c=1.5$. We also assume that $b=20$. Then for different values of $\theta$ we have

\begin{itemize}
\item $\theta = -1$
\begin{eqnarray}
\chi(u,b) = 1 + 0.0186 e^{-2.2207u} - 0.7223 e^{-0.2687u},\nonumber
\end{eqnarray}
\item $\theta = -0.5$
\begin{eqnarray}
\chi(u,b) = 1 + 0.0105 e^{-2.1148u} - 0.6970 e^{-0.2976 u}\nonumber,
\end{eqnarray}
\item $\theta=0.5$
\begin{eqnarray}
\chi(u,b) = 1 - 0.0140 e^{-1.8736 u} - 0.6314 e^{-0.3788 u}\nonumber,
\end{eqnarray}
\item $\theta=1$
\begin{eqnarray}
\chi(u,b) = 1 - 0.0335 e^{-1.7305 u} - 0.5866 e^{-0.4392 u}\nonumber 
\end{eqnarray}
\end{itemize}

\end{example}

\section{Discussion}
In this paper we derived existing results in the literature by implementing simple methods used in the classical risk model. We considered the situation that inter-claim times and claim amounts are dependent. We assumed that such dependence structure can be explained through FGM copula. Unfortunately, it is difficult to find such expressions for other copula functions. Further, the FGM copula cannot demonstrate a strong relationship between random variables as for example, Clayton copula. Incorporation of other copula functions to risk models and the approximation of the ruin probability can be the subject of other studies.

\end{document}